\documentclass[12pt]{amsart}
\usepackage{amscd}
\usepackage{amsmath}
\usepackage{amssymb}

\usepackage[all]{xy}
\def\frk{\frak}               

\def\Phi{{\frk n}}
\def\Phi{{\frk N}}
\def\opn#1#2{\def#1{\operatorname{#2}}} 
\opn\chara{char} \opn\length{\ell} \opn\pd{pd} \opn\rk{rk}
\opn\projdim{proj\,dim} \opn\injdim{inj\,dim} \opn\rank{rank}
\opn\depth{depth} \opn\sdepth{sdepth} \opn\fdepth{fdepth}
\opn\grade{grade} \opn\height{height} \opn\embdim{emb\,dim}
\opn\codim{codim}  \opn\min{min} \opn\max{max}

\opn\Tr{Tr} \opn\bigrank{big\,rank}
\opn\superheight{superheight}\opn\lcm{lcm}
\opn\trdeg{tr\,deg}
\opn\reg{reg} \opn\lreg{lreg} \opn\ini{in} \opn\lpd{lpd}
\opn\size{size}
\opn\div{div} \opn\Div{Div} \opn\cl{cl} \opn\Cl{Cl}
\opn\Spec{Spec} \opn\Supp{Supp} \opn\supp{supp} \opn\Sing{Sing}
\opn\Ass{Ass} \opn\Min{Min}
\opn\Ann{Ann} \opn\Rad{Rad} \opn\Soc{Soc}
\opn\Im{Im} \opn\Ker{Ker} \opn\Coker{Coker} \opn\Am{Am}
\opn\Hom{Hom} \opn\Tor{Tor} \opn\Ext{Ext} \opn\End{End}
\opn\Aut{Aut} \opn\id{id}  \opn\deg{deg}

\opn\nat{nat}
\opn\pff{pf}
\opn\Pf{Pf} \opn\GL{GL} \opn\SL{SL} \opn\mod{mod} \opn\ord{ord}
\opn\Gin{Gin} \opn\Hilb{Hilb}
\opn\aff{aff} \opn\con{conv} \opn\relint{relint} \opn\st{st}
\opn\lk{lk} \opn\cn{cn} \opn\core{core} \opn\vol{vol}
\opn\link{link} \opn\star{star}
\opn\gr{gr}

\def\pot#1#2{#1[\kern-0.28ex[#2]\kern-0.28ex]}

\opn\dirlim{\underrightarrow{\lim}}
\opn\inivlim{\underleftarrow{\lim}}
\let\to=\rightarrow

\def\Implies{\ifmmode\Longrightarrow \else
        \unskip${}\Longrightarrow{}$\ignorespaces\fi}
\def\implies{\ifmmode\Rightarrow \else
        \unskip${}\Rightarrow{}$\ignorespaces\fi}
\def\iff{\ifmmode\Longleftrightarrow \else
        \unskip${}\Longleftrightarrow{}$\ignorespaces\fi}

\let\:=\colon
\newtheorem{Theorem}{Theorem}[section]

\newtheorem{Remark}[Theorem]{Remark}

\newtheorem{Example}[Theorem]{Example}

\let\epsilon\varepsilon
\let\phi=\varphi
\let\kappa=\varkappa
\def\qed{\ifhmode\textqed\fi
      \ifmmode\ifinner\quad\qedsymbol\else\dispqed\fi\fi}
\def\textqed{\unskip\nobreak\penalty50
       \hskip2em\hbox{}\nobreak\hfil\qedsymbol
       \parfillskip=0pt \finalhyphendemerits=0}
\def\dispqed{\rlap{\qquad\qedsymbol}}

\opn\dis{dis}
\def\pnt{{\raise0.5mm\hbox{\large\bf.}}}

\opn\Lex{Lex}

\begin{document}
\title{\bf Depth of  some special  monomial  ideals}

\author{ Dorin Popescu and Andrei Zarojanu}

\thanks{The  support from  grant ID-PCE-2011-1023 of Romanian Ministry of Education, Research and Innovation is gratefully acknowledged.}

\address{Dorin Popescu,  "Simion Stoilow" Institute of Mathematics of Romanian Academy, Research unit 5,
 P.O.Box 1-764, Bucharest 014700, Romania}
\email{dorin.popescu@imar.ro}
\address{Andrei Zarojanu,  Faculty of Mathematics and Computer Sciences, University
of Bucharest, Str. Academiei 14, Bucharest, Romania, and}
\address{"Simion Stoilow" Institute of Mathematics of Romanian Academy, Research group of the project  ID-PCE-2011-1023,
 P.O.Box 1-764, Bucharest 014700, Romania}
\email{andrei\_zarojanu@yahoo.com}

\maketitle
\maketitle
\begin{abstract} Let $I\supsetneq J$ be  two  squarefree monomial ideals of a polynomial algebra over a field. Suppose that $I$ is generated by one squarefree monomial of degree $ d>0$, and other squarefree monomials  of degrees $\geq d+1$. If  the Stanley depth of $I/J$ is $\leq d+1$ then almost always  the usual depth of $I/J$ is $\leq d+1$ too.

\end{abstract}

\section*{Introduction}

 Let $S=K[x_1,\ldots,x_n]$ be the polynomial algebra in $n$ variables over a field $K$ and   $I\supsetneq J$  two   squarefree monomial ideals of $S$. Suppose that $I$ is generated by squarefree monomials of degrees $\geq d$   for some positive integer $d$. Modulo a multigraded isomorphism we may assume either that $J=0$, or $J$ is generated in degrees $\geq d+1$.  Then $\depth_S I/J\geq d$ (see  \cite[Proposition 3.1]{HVZ}, \cite[Lemma 1.1]{P}) and  upper bounds  are given by Stanley's Conjecture if it  holds. Finding such upper bounds is the subject of several papers  \cite{P2}, \cite{P}, \cite{P1}, \cite{Sh}, \cite{PZ}.  We remind  below the notion  of Stanley depth.

   Let $P_{I\setminus J}$  be the poset of all squarefree monomials of $I\setminus J$  with the order given by the divisibility. Let ${\mathcal P}$ be a partition of  $P_{I\setminus J}$ in intervals $[u,v]=\{w\in  P_{I\setminus J}: u|w, w|v\}$, let us say   $P_{I\setminus J}=\cup_i [u_i,v_i]$, the union being disjoint.
Define $\sdepth {\mathcal P}=\min_i\deg v_i$ and  the  {\em Stanley depth} of $I/J$ given by $\sdepth_SI/J=\max_{\mathcal P} \sdepth {\mathcal P}$, where ${\mathcal P}$ runs over the set of all partitions of $P_{I\setminus J}$ (see  \cite{HVZ}, \cite{S}).
 Stanley's Conjecture says that  $\sdepth_S I/J\geq \depth_S I/J$.

Let $r$ be the number of  squarefree monomials of degree $d$ of $I$ and $B$ (resp. $C$) be the set of  squarefree monomials of degrees $d+1$  (resp. $d+2$) of $I\setminus J$.  Set $s=|B|$, $q=|C|$. If either $s>r+q$, or $r>q$, or $s<2r$ then $\sdepth_SI/J\leq d+1$ and if  Stanley's Conjecture holds then any of these numerical conditions would imply $\depth_SI/J\leq d+1$, independently of the characteristic of $K$. In particular this was proved  directly in \cite{P1} and \cite{Sh}.

Suppose that $r=1$. If $d=1$ we showed in \cite[Theorem 1.10]{PZ} that almost always  Stanley's Conjecture holds. It is the purpose of this note to complete the result  for $d\geq 1$  in the next form.

{\bf Theorem} {\em Suppose that $I \subset S$ is minimally generated by a squarefree monomial $f$ of degree $d$,  and a  set $E$ of squarefree monomials of degrees $\geq d+1$. Assume that $s\not =q+1$ and  $\sdepth_S I/J=d+1$.
Then $\depth_S I/J \leq d+1$.}

\section{Proof of the Theorem }

We may assume  that  $s<q+1$ because if $s>q+1 $ then by \cite{P1} we get $\depth_S I/J\leq d+1$. Also we may suppose that   $C\subset (f,B)$ by \cite[Lemma 1.6]{PZ}.
  Induct on $|E|$. Suppose that $E=\emptyset$. If $C=\emptyset$, then   $\depth_SI/J\leq d+1$ by \cite[Lemma 1.5]{PZ}. Otherwise, let $c=fx_{n-1}x_n\in C$ and $I'=(B\setminus \{fx_{n-1},fx_{n}\})$. In the exact sequence
  $$0\to I'/J\cap I'\to I/J\to I/J+I'\to 0$$
 the last term has sdepth $d+2$ since $c\not \in I'+J$ and so the first one has sdepth $\leq d+1$ by \cite[Lemma 2.2]{R} and even depth $\leq d+1$ by \cite[Theorem 4.3]{P}. Then the Depth Lemma gives $\depth_SI/J\leq d+1$.

   Set $I_n=(B \setminus \{fx_n\}),J_n=I_n\cap J$. In the following exact sequence
   $$0\to I_n/J_n\to I/J\to I/(I_n+J)\to 0$$
   the last term has sdepth depth $d+1$ since $[f,fx_n]$ is the whole poset of $(f)/(f)\cap (I_n+J)$ and $x_n \not \in ((J+I_n):f)$. If the first term has sdepth $=d+1$ then by [\cite{P}, Theorem 4.3.] we get depth $=d+1$ and applying Depth Lemma the conclusion follows. So we can assume that there exists  a partition ${\mathcal P}_n$  of $I_n/J_n$ with sdepth $d+2$. We may suppose that all intervals of ${\mathcal P}_{b_i}$ (as well as of other partitions which we will use)  starting with a monomial $v$ of degree $\geq d+2$ have the form $[v,v]$.
 In ${\mathcal P}_n$ we can't have the interval $[c,c]$, $c=fx_{n-1}x_n$, or the interval $[fx_{n-1},c]$ because otherwise we can switch it with $[f,c]$ and get a partition of $I/J$ with sdepth $d+2$. Thus we have in ${\mathcal P}_n$ the interval $[b_1,c]$, $b_1 \in E$. Switching the interval $[b_1,c]$ with the interval $[fx_n,c]$ we get a partition $P_{B_{b_1}}$ for $(B_{b_1})/J_{b_1}$ where $B_{b_1} = B \setminus \{b_1\}$ and $J_{b_1}=(B_{b_1}) \cap J$.

  In ${\mathcal P}_{B_{b_1}}$ we have an interval $[\bar{c},\bar{c}]$ because $s < q+1$. Thus there exists $b_2 \in E$ such that $\bar{c} \in (b_2)$. Note that $ \bar{c} \notin (b_1)$ because otherwise we may replace in $ {\mathcal P}_{B_{b_1}}$ the interval $[\bar{c},\bar{c}]$ with the interval $[b_1,\bar{c}]$ and get a partition ${\mathcal P}_B$ for $(B) / (B) \cap J$ with sdepth $=d+2$. This leads to a contradiction because we may change in $P_B$ two intervals like $[fx_i,fx_ix_j]$, $[fx_j,c']$ for some $c'$ with $[f,fx_ix_j]$, $[c',c']$ and get a partition ${\mathcal P}$ for $I/J$ with sdepth $d+2$.

   Let $I_{b_2}=(f,E\setminus \{b_2\})$, $J_{b_2}=I_{b_2}\cap J$. In the following exact sequence
  $$0\to I_{b_2}/J_{b_2}\to I/J\to I/(I_{b_2}+J)\to 0$$
the last term has depth $\geq d+1$ because it is isomorphic with $(b_2)/(b_2)\cap (I_{b_2}+J)$. If $\sdepth I_{b_2} / J_{b_2} \leq d+1$ then the first term has depth $\leq d+1$ by the induction hypothesis, so by the Depth Lemma we get $\depth I/J \leq d+1$.

  Now assume that $ \sdepth I_{b_2} / J_{b_2} \geq d+2$ and let $P_{b_2}$ be a partition on $I_{b_2}/J_{b_2}$ with
sdepth $d+2$. In ${\mathcal P}_{b_2}$ we have the interval $[f,fx_ix_j],i,j \in [n] \setminus \mbox{supp}f$. We have in ${\mathcal P}_{b_2}$  for all $b \in B \setminus \{b_2,fx_i,fx_j\}$ an interval $[b,c_b]$. We define $h_2:(B\setminus \{b_2\})\to C$ by $b\to c_b$ and $h_2 (fx_i)=h_2(fx_j)=fx_ix_j$ and let $g_2:\Im h_2 \to  (B\setminus \{b_2\})$ defined by $c_b\to b$, $g_2 (fx_ix_j)=fx_i$. Similarly we define $h_1,g_1$ for ${\mathcal P}_{B_{b_1}}$, that is $h_1$ is given by $b'\to c'$ if ${\mathcal P}_{B_{b_1}}$ has the interval $[b',c']$.

We want to show that we can build a partition ${\mathcal P}$ with sdepth $=d+2$ for $I/J$. Consider $a_0 = b_1$ and $ c_{i-1} = h_2(a_{i-1}), a_i = g_1(c_{i-1}), i > 0$.
The construction stops at step $e$ if

1) $a_e=b_2$,

2) $c_e \notin \Im h_1$,

3)$a_e=fx_{j}$ after $a_u=fx_{i}, u < e$ already appeared. Note that here we have fixed $a_u = fx_i$.

In the first case we set $c_e=\bar{c}$ and we see that $h_1$ gives a bijection between $\{a_1,\ldots,a_e\}$ and $\{c_0,\ldots,c_{e-1}\}$. But $h_1$ also gives a bijection between $B\setminus \{b_1,a_1,\ldots,a_e\}$ and $C\setminus \{\bar{c},c_0,\ldots,c_{e-1}\}$. Then the intervals $[a_p,c_p]$, $0\leq p\leq e$ and the intervals $[g_1({\tilde c}),{\tilde c}]$, ${\tilde c}\in C\setminus \{\bar{c},c_0,\ldots,c_{e-1}\}$ and some other intervals starting with monomials of degree $\geq d+2$ give a partition  ${\mathcal P}_B$ of $P_{B/B \cap J}$ with sdepth $\geq d+2$. As before this is a contradiction with sdepth $I/J = d+1$.

In the second case, as above we see that the intervals $[a_p,c_p]$, $0\leq p\leq e$ and the intervals $[g_1({\tilde c}),{\tilde c}]$, ${\tilde c}\in C\setminus \{\bar{c},c_0,\ldots,c_{e-1}\}$ and some other intervals starting with monomials of degree $\geq d+2$ give a partition  ${\mathcal P}_B$ of $P_{B/B \cap J}$ with sdepth $\geq d+2$. Contradiction.

In the last case we see as usual, that $h_1$ gives a bijection between $\{a_1,\ldots,a_e\}$ and $\{c_0,\ldots,c_{e-1}\}$. But $h_1$ also gives a bijection between $B\setminus \{b_1,a_1,\ldots,a_e\}$ and $C\setminus \{c_0,\ldots,c_{e-1}\}$. Then the intervals $[a_p,c_p]$, $0\leq p\leq e-1, p \neq u$ and the intervals $[f,c_u],[g_1({\tilde c}),{\tilde c}]$, ${\tilde c}\in C\setminus \{c_0,\ldots,c_{e-1}\}$ and some other intervals starting with monomials of degree $\geq d+2$ give the partition  ${\mathcal P}$ of $P_{I/J}$ with sdepth $\geq d+2$. Contradiction.

   \vspace{30 pt}

\begin{Example} {\em
Let $n=5$, $I = (x_1x_2,x_3x_4x_5)$ and $J=(x_1x_2x_3x_5,x_1x_2x_4x_5)$. We see that we have $\sdepth I/J = d+1=3$ and $B = \{ x_1x_2x_3,x_1x_2x_4,x_1x_2x_5,x_3x_4x_5 \}$ and $C = \{x_1x_2x_3x_4, x_1x_3x_4x_5, x_2x_3x_4x_5 \}$ so we are in the case $s=q+1$.
We can get $\depth I/J \leq 3$   by
 using \cite[Lemma 1.5]{PZ}  for $u=x_1x_2x_5$. }
\end{Example}

\begin{Remark} {\em If in the above example change just one monomial from the generators of $J$, namely take $J=(x_1x_2x_4x_5,x_2x_3x_4x_5)$ then we have $\sdepth_SI/J=4$ because the partition induced by the intervals $[x_1x_2,x_1x_2x_3x_4]$, $[x_3x_4x_5,x_1x_3x_4x_5]$, $[x_1x_2x_5,x_1x_2x_3x_5]$ has sdepth $d+2=4$. Also we have $\depth_SI/J=4$. }
\end{Remark}
A question is hinted by the following example.
\begin{Example} {\em
Let $n=5$, $I = (x_1, x_2x_3,x_2x_4,x_2x_5,x_3x_4)$ and $J$ the ideal generated by all squarefree monomials of $I$ of degrees $4$. Then $E=\{ x_2x_3,x_2x_4,x_2x_5,x_3x_4\}$, $f=x_1$,  $B=\{x_1x_2,x_1x_3,x_1x_4,x_1x_5,E\}$, $C=\{x_1x_2x_3,x_1x_2x_4,x_1x_2x_5,x_1x_3x_4,$\\ $x_2x_3x_4,
x_2x_3x_5,x_2x_4x_5\}$. Thus $s=8=q+1$. We see that  $\sdepth_SI/J=d+2=3$ but $\depth_SI/J=d+1=2$. Note that here $$C\subset (\cup_{a\in E} C\cap (f)\cap (a))\cup (\cup_{a,a'\in E, a\not =a'} C\cap (a)\cap (a'),$$
a condition which might imply always $\depth_SI/J\leq d+1$, the inequality being not true for sdepth. }
\end{Example}

\end{document}